\documentclass{article}
\usepackage{pdflscape}
\usepackage{
 amsmath,
 amsxtra,
 amsthm,
 amssymb,
 etex,
 mathrsfs,
 }
\usepackage[all]{xy}

\newtheorem{theorem}{Theorem}[section]
\newtheorem{lemma}[theorem]{Lemma}

\newtheorem{proposition}[theorem]{Proposition}
\newtheorem{corollary}[theorem]{Corollary}

\theoremstyle{definition}
\newtheorem{defn}[theorem]{Definition}
\newtheorem{remark}[theorem]{Remark}

\newcommand{\bd}{\begin{defn}}
\newcommand{\ed}{\end{defn}}
\newcommand{\bl}{\begin{lemma}}
\newcommand{\el}{\end{lemma}}
\newcommand{\bp}{\begin{proposition}}
\newcommand{\ep}{\end{proposition}}
\newcommand{\bt}{\begin{theorem}}
\newcommand{\et}{\end{theorem}}
\newcommand{\bc}{\begin{corollary}}
\newcommand{\ec}{\end{corollary}}
\newcommand{\br}{\begin{remark}}
\newcommand{\er}{\end{remark}}
\newcommand{\bpf}{\begin{proof}}
\newcommand{\epf}{\end{proof}}

\newcommand{\Z}{\mathbb{Z}}
\newcommand{\Q}{\mathbb{Q}}

\newcommand{\Op}{\mathcal{O}}
\newcommand{\al}{\alpha}
\newcommand{\be}{\beta}

\newcommand{\ze}{\zeta}

\DeclareMathOperator{\Gal}{Gal}

\newcommand{\lra}{\longrightarrow}

\numberwithin{equation}{section}

\begin{document}
\title{On the $p$-divisibility of even $K$-groups of the ring of integers of a cyclotomic field}
 \author{
  Meng Fai Lim\footnote{School of Mathematics and Statistics, Key Laboratory of Nonlinear Analysis and Applications (Ministry of Education),
Central China Normal University, Wuhan, 430079, P.R.China.
 E-mail: \texttt{limmf@ccnu.edu.cn}} }
\date{}
\maketitle

\begin{abstract} \footnotesize
\noindent
Let $k$ be a given positive odd integer and $p$ an odd prime. In this paper, we shall give a sufficient condition when a prime $p$ divides the order of the groups $K_{2k}(\Z[\ze_m+\ze_m^{-1}])$ and $K_{2k}(\Z[\ze_m])$, where $\ze_m$ is a primitive $m$th root of unity. When $F$ is a $p$-extension contained in $\Q(\ze_\ell)$ for some prime $\ell$, we also establish a necessary and sufficient condition for the order of $K_{2(p-2)}(\Op_F)$ to be divisible by $p$, which generalizes a previous result of Browkin.


\end{abstract}

\section{Introduction}

Let $F$ be a number field, and $\Op_F$ its ring of integers. Thanks to the results of Quillen \cite{Qui73b}, Garland \cite{Gar} and Borel \cite{Bo}, the even $K$-group $K_{2k}(\Op_F)$ is known to be finite for every $k\geq 1$. A conjecture of Birch and Tate \cite{Bir, Tate} predicted that the exact size of $K_2(\Op_F)$ is related to the value of the Dedekind zeta function of $F$ at $s=-1$. This conjecture was subsequently generalized by Lichtenbaum \cite{Lic72, Lic73} to higher $K$-groups. Thanks to the effort of many, this conjecture of Birch, Tate and Lichtenbaum is now known to be valid for a totally real abelian field (see Theorem \ref{Lich} for a description of this).

Building on these known formulas, many authors have studied and calculated the order of $K_2$-groups of various totally real abelian fields (see \cite{BT, Bir, Brow00, Brow05, BrowGang, Zh}, where this list is probably far from being exhaustive). Recently, the author together with Qin \cite{LimQin} has performed numerical computations and obtained the precise orders of higher even $K$-groups for various totally real abelian fields, again building on known cases of Lichtenbaum's conjecture.

In this paper, we shall take a different perspective of study from the above mentioned works. Namely, we like to ask whether one can obtain any theoretical information on the $K$-groups via known cases of conjectures of Birch-Tate and Lichtenbaum. Our modest attempt towards this to establish certain concrete divisibility properties of the $K_{2k}(\Op_F)$ for number fields contained in a cyclotomic field. We should mention that the divisibility properties of $K_2$ has been studied by many (for instances, see \cite{Brow85, Brow05, Ku, Lu, Wu, Zh}).

To state our result, we need to introduce certain notation that will be adhered to throughout the article. Let $p$ be an odd prime. For a given positive integer $N$, $\ze_N$ will always denote a primitive $N$th root of unity. For every integer $n$, we then write $v_p(n)$ for the largest integer such that $p^{v_p(n)}$ divides $n$.
In the event that $v_p(n)=1$, we sometimes write $p\parallel n$. For a ring $R$ with identity, we shall write $K_n(R)$ for the
algebraic $K$-groups of $R$ in the sense of Quillen \cite{Qui73b} (also see \cite{KolTrieste, WeiKbook}).

We can now state the main theorem of the paper.

\bt \label{mainintro} Suppose that $m$ is an integer $>1$ such that at least one of its prime divisor satisfies $\equiv 1$ (mod $p$). Set $d^+=|\Q(\ze_p,\ze_m+\ze_m^{-1}):\Q(\ze_m+\ze_m^{-1})|$ and $d=|\Q(\ze_p,\ze_m):\Q(\ze_m)|$.

Then for every pair of integers $(r,k)$ such that $r \geq 1$ and $k$ a positive odd integer $<p-2$, the orders of the groups $K_{2k+2d^+r}(\Z[\ze_m+\ze_m^{-1}])$ and $K_{2k+2dr}(\Z[\ze_m])$ are both divisible by $p$.

In the event that the integer $m$ has at least one prime divisor satisfying $\equiv 1$ (mod $p^2$), we also have $p$ dividing both the orders of $K_{2(p-2)+2d^+r}(\Z[\ze_m+\ze_m^{-1}])$ and $K_{2(p-2)+2dr}(\Z[\ze_m])$.
\et

Note that in the theorem, one has
\[ d^+=\left\{
         \begin{array}{ll}
           p-1, & \hbox{if $p \nmid m$;} \\
           2, & \hbox{if $p \mid m$,}
         \end{array}
       \right. \quad\mathrm{and}\quad  d=\left\{
         \begin{array}{ll}
           p-1, & \hbox{if $p \nmid m$;} \\
           1, & \hbox{if $p \mid m$.}
         \end{array}
       \right.
\]

In the case when $k\leq p-2$, we can even obtain a lower bound on the power of $p$ dividing the order of the $K_{2k}$-groups.

\bt[Theorem \ref{main}] Let $p$ be an odd prime, $k$ a positive odd integer and $m$ an integer $>1$. Define
\[s_j = \#\big\{\ell \mbox{ is a prime divisor of } m~\big|~ v_p(\ell-1)=j\big\} \]
and set $\theta:=\max\{j~|~s_j\neq 0\}$. Suppose that $p\geq k+2$.
Then the power of $p$ dividing the orders of the groups $K_{2k}(\Z[\ze_m+\ze_m^{-1}])$ is at least
\[\frac{p^{s_1+s_2\cdots +s_\theta}-p^{\delta s_1}+\delta}{p-1} + \sum_{j=2}^{\theta}\frac{p^{s_j+\cdots +s_\theta}-1}{p-1}p^{s_1+2s_2+\cdots + (j-1)s_{j-1} + (j-1)(s_j+\cdots s_\theta) -(j-1)},\]
where
\[ \delta =\left\{
     \begin{array}{ll}
       0, & \hbox{if $p>k+2$;} \\
       1, & \hbox{if $p=k+2$.}
     \end{array}
   \right.\]
A similar divisibility assertion also holds for $K_{2k}(\Z[\ze_m])$.
\et

In particular, by the preceding theorem, we see that \[ 7^{1+7+7^2} \mid \#K_{6}(\Z[\zeta_{29\cdot 43\cdot 71}])\]
(see Remark \ref{remark div} for more discussion on this).

For the remainder of the introductory section, we describe the outline of the proof leaving details to the body of the paper. For simplicity, we shall first assume in this introduction that $k=1$ and $m=\ell$, where $\ell$ is a prime  $\equiv 1$ (mod $p$) (see the body of the paper for the general statements).
The first step is to apply certain descent arguments (see Subsection \ref{descent results} and proof of Theorem \ref{main}) to reduce the divisibility assertion of our theorem to a divisibility assertion over a $p$-extension $F$ contained in $\Q(\ze_\ell)$. As our prime $p$ is odd, such an extension $F$ is necessarily totally real and abelian, and so we may apply the known case of Birch, Tate and Lichtenbaum's conjecture here. Via Artin formalism and followed by certain intermediate argument, one obtains
\[\# K_{2}(\Op_F) \sim_p \prod_{\substack{\chi\in \widehat{\Gal(F/\Q)}\\ \chi\neq\chi_0}}B_{2,\chi},\]
where $\chi_0$ denotes the trivial character. Here, if $\al,\be \in \overline{\Q}$, we write $\al\sim_p\be$ to mean that $\al/\be$ is a $p$-adic unit. We are therefore reduced to proving $p$-divisibility of the generalised Bernoulli number $B_{2,\chi}$ for a non-trivial character $\chi$ of $\Gal(F/\Q)$. Recall that by \cite[Exercise 4.2(a)]{Wash}, one has
\[ B_{2,\chi} = -\ell^{-1}\sum_{a=1}^{\ell-1}\chi(a)a^2.  \]
To proceed, we utilize an idea going back to Browkin \cite{Brow05} and Wu \cite{Wu}. Namely, since $\Gal(F/\Q)$ is a cyclic $p$-group, the character $\chi$ will take values in $\Z[\ze_{p^N}]$ for some big enough integer $N$. Hence we have
\[ B_{2,\chi} \equiv -\frac{~1~}{~\ell~}\sum_{a=1}^{\ell-1}a^2 = -\frac{(\ell-1)(2\ell-1)}{6}  \quad\mathrm{mod}~(1-\ze_{p^N})\Z[\ze_{p^N}].\]
Here the first congruence makes sense since $\ell$ is coprime to $p$ (and hence $(1-\ze_{p^N})\Z[\ze_{p^N}]$).
Plainly, the final expression $\equiv 0$ if either $p>3$ or $p=3$ with $p^2\mid (\ell-1)$.
Consequently, we have $p$ dividing $\# K_{2}(\Op_F)$ whenever either $p>3$ or $p=3$ with $p^2\mid (\ell-1)$. In fact, one can even obtain a lower bound on the power of $p$ dividing $\# K_{2}(\Op_F)$ by a finer argument (see Theorem \ref{main}).

We now say a bit for a general odd $k>1$ and $m$. The same idea as above carries over except that the process is slightly more involved. For this, we need to be able to control $p$-contribution coming from the denominators of Bernoulli numbers. This part of the argument requires us to work under the assumption that $k\leq p-2$, where one has take slightly more care on the ``boundary condition'' $k=p+2$. For instance, as seen above, when $p=3$ and $k=1$, the term $3$ appears in the denominator term after applying the power sum formulas. Finally, we need to analyse the $p$-divisibility of $B_{k+1,\chi}$ which led us to consider sums of the form \[\sum_{\substack{a=1 \\ \gcd(a,\ell_1\cdots \ell_s)=1}}^{\ell_1\cdots \ell_s-1} a^{k+1-i},\] where $\ell_1,...,\ell_s$ are distinct primes satisfying $\equiv 1$ (mod $p$). Plainly, the power sum formula cannot apply to this kind of sum on the nose. But thankfully, a simple manipulation trick (see (\ref{B=powersum2})) allows us to rewrite this sum as a summation of appropriate power sums. From here, we apply certain properties of the general power sum formula (see Lemma \ref{power sum divi}) to obtain $p$-divisibility. Again, throughout the argument process, we need to take care on the possible $p$-contribution from the denominator in the power sum formula. To handle this latter issue, we require to call upon the work of Kellner and Sondow \cite{KS}, and we refer readers to Section \ref{Main results section} for the details.

Finally, if $k>p+2$, it does not seem easy to determine the $p$-divisibility of $B_{k+1,\chi}$ due to possibly multiple appearances of large $p$-power in denominators. For this, we circumvent this by appealing to a periodicity result on the $p$-rank of $K$-groups (see Proposition \ref{periodicity of p-rank}) which allows us to obtain $p$-divisibility of the $K_{2k}$-groups for appropriate $k>p+2$.

We finally specialize to cyclic degree $p$-extensions contained in $\Q(\ze_\ell)$ for some prime $\ell$ at the ``boundary condition'' $k=p+2$, where we establish a generalization of Browkin's result (see Theorem \ref{BrowkinK}). In fact, we can give a proportion when the divisibility is valid.

\bc[Corollary \ref{BrowkinKcor}]
Let $p$ be an odd prime. For every prime $\ell\equiv 1$ $($mod $p)$, we write $F_\ell$ for the  cyclic extension of $\Q$ of degree $p$ contained in $\Q(\ze_\ell)$. Then one has
\[\lim_{x\to\infty}\frac{\#\{\ell~prime \leq x~|~\mbox{$p$ divides $\#K_{2(p-2)}(\Op_{F_\ell})$ and $\ell\equiv 1$ $($mod $p)$}\}}{\#\{\ell~prime \leq x~|~\mbox{$\ell\equiv 1$ $($mod $p)$}\}} =\frac{~1~}{~p~}.\]
\ec

Note that $p$ does not divide $\#K_{2(p-2)}(\Op_{F_\ell})$ if and only if the number field $F_\ell$ is $(p, p-1)$-regular in the sense of Assim \cite{Ass}. The latter is equivalent to saying that $F$ is $p$-rational in the sense of Movahhedi \cite{Mov, MN} (see \cite[Sections
4.2 and 6]{AssMov}; we thank the anonymous referee for pointing out this observation). A recent work of Assim-Movahhedi \cite{AssMov} has provided a complete list of $p$-extensions $F$ of $\Q$ with $K_{2k}(\Op_F)[p^\infty]=0$ for every $k$ (we once again thank the anonymous referee for bringing this paper to our attention). Therefore, our Corollary \ref{BrowkinKcor} will also follow as a corollary of the work of Assim-Movahhedi, whose approach is different from the one we have adopted here.

We like to further mention a perspective of our study. The question of the $p$-divisibility of $B_{k+1,\chi}$ can be thought as a question of whether the Dirichlet $L$-value $L(\chi, -k)$ vanishes modulo $p$. The recent work of Burungale and Sun \cite{BS} studied this latter problem, although the characters they considered have conductors $f$ satisfying $p\nmid \phi(f)$ which is in contrast to what we do here. Furthermore, in our work, we study the generalized Bernoulli number $B_{k+1,\chi}$ directly, whereas their investigation revolves around understanding the $L$-values. Therefore, our work may be thought of a complement to their in terms of both the objects of study and investigation approaches.

We end the introductory section giving an outline of the paper. In Section \ref{Preliminaries}, we collect certain facts about the power sum denominator and  Bernoulli numbers with emphasis on the possible $p$-contribution in the power sum denominator and the denominator of Bernoulli numbers. Section \ref{descent results} is where we recall two descent results and a periodicity result of $K$-groups. We then recall the important known case of Lichtenbaum's conjecture in Section \ref{Lich sec} which will be an important tool for the eventual proof of our main theorem.
In Section \ref{Main results section}, we will prove our main theorem. We also give a lower bound on the power of $p$ dividing the $K_{2k}$-groups when $k\leq p-2$.
In Section \ref{BrowkinSec}, we revisit a result of Browkin and generalize this said result to higher $K$-groups. Finally, in Section \ref{example sec}, we provide some numerical examples to illustrate our results.

\subsection*{Acknowledgement}

The author likes to thank Ashay Burungale, Chao Qin and Jun Wang for many insightful discussions during the preparation of the paper. He would also like to express his gratitude to Ashay Burungale for providing thorough responses to the numerous questions he had regarding his paper \cite{BS}. This research is supported by the
National Natural Science Foundation of China under Grant No.\ 11771164. Finally, the author wishes to express sincere gratitude to the anonymous referee for providing valuable comments and suggestions that have significantly enhanced the clarity of the manuscript. Additionally, the author likes to
thank the anonymous referee for making him aware of the paper \cite{AssMov}.

\section{Power-sum denominator and Bernoulli number} \label{Preliminaries}
We begin collecting certain facts about the power sum denominator and and Bernoulli number that will be required for the discussion in this paper. As a start, one recalls that the Bernoulli numbers $B_n$ are defined by
\[ \frac{t}{e^t-1} = \sum_{n=0}^\infty B_n\frac{t^n}{n!}.\]

\bl \label{Bernoulli denom}
Let $n$ be a positive even integer. Then the denominator of the Bernoulli number $B_n$ is given by
\[ \prod_{\substack{p~\mathrm{prime}\\ (p-1)\mid n}} p .\]
In particular, the denominator of the Bernoulli number $B_n$ is squarefree with every prime divisor $\leq n+1$.
\el

\bpf
This is an immediate consequence of the von Staudt-Clausen Theorem (see \cite[Theorem 5.10]{Wash}).
\epf

We also recall the $n$th Bernoulli polynomial $B_n(x)$ which are defined by
\[ \frac{te^{xt}}{e^t-1} = \sum_{n=0}^\infty B_n(x)\frac{t^n}{n!}.\]
The relation between the Bernoulli numbers and the Bernoulli polynomials is encoded in the following formula
\begin{equation}\label{Bern poly}
B_n(x) = \sum_{i=0}^n \binom{n}{i}B_i x^{n-i}.
\end{equation}

\bd
Set $S_n(x) : = \displaystyle\frac{1}{n+1}\big(B_{n+1}(x) - B_{n+1} \big)$.
\ed

Clearly, $S_n(x)$ is a polynomial with rational coefficients. Furthermore, for every positive integer $m$, we have
\begin{equation} \label{power sum and Bernoulli}
S_n(m) = \frac{1}{n+1}\big(B_{n+1}(m) - B_{n+1} \big) = \sum_{a=0}^{m-1} a^n.
\end{equation}

\bd For $n \geq 1$, the $n$th power-sum denominator is defined to be the smallest positive
integer $d_n$ such that $d_n S_n(x)$ is a polynomial in $x$ with integer
coefficients. We also set $d_0 =1$.

Some of the values of $d_n$ are listed in \cite[Sequence A064538]{Sl}.
\ed

\br
In view of (\ref{power sum and Bernoulli}), $d_n$ can be thought as the smallest positive
integer $d_n$ such that the following formal expression
\[\mbox{``} d_n (1^{x-1}+\cdots+ a^{x-1}) \mbox{''} \]
 is a polynomial in $x$ with integer
coefficients. This is the definition given in \cite[Definition 1]{KS}.

\er

We now record certain divisibility properties of these power-sum denominators that will be required for subsequent discussion. (For more detailed account on this topic, we refer readers to the article of Kellner and Sondow \cite{KS} and the reference therein.) Now, if $p$ is a prime and $N$ is an integer, we write $p\parallel N$ to mean $p\mid N$ but $p^2\nmid N$.

\bp \label{powersumdivisor}
For every $n\geq 1$, the following statements are valid.
\begin{enumerate}
  \item[$(a)$]  Every prime divisor of $d_n$ is $\leq n+1$.
  \item[$(b)$]  The integer $d_n$ is divisible by $n+1$ and $\displaystyle\frac{d_n}{n+1}$ is square-free. Furthermore, every prime divisor of $\displaystyle\frac{d_n}{n+1}$ is less than or equal to $M_n$, where
      \[ M_n=\left\{
              \begin{array}{ll}
               \displaystyle \frac{n+2}{2}, & \hbox{if $n$ is even;} \\
               \\
                \displaystyle\frac{n+2}{3}, & \hbox{if $n$ is odd.}
              \end{array}
            \right.
      \]
  \item[$(c)$]  If $n+1$ is a prime, then $(n+1) \parallel d_n$.
\end{enumerate}
\ep

\bpf
Assertions (a) and (b) follow from \cite[Theorem 1]{KS}. Finally, in view that $M_n<n+1$, assertion (c) is an immediate consequence of assertion (b).
\epf

We now like to understand the polynomial $S_n(x)$ more. In view of (\ref{power sum and Bernoulli}), it is then natural to start looking at the power sum $\sum_{a=0}^{m-1} a^n$. To get a feel, we begin writing the power sum formulas for a few values of $n$'s.

\begin{eqnarray*}
  \sum_{a=0}^{m-1} a  &=& \frac{1}{2}(m-1)m \\
  \sum_{a=0}^{m-1} a^2  &=&  \frac{1}{6}(m-1)m(2m-1) \\
  \sum_{a=0}^{m-1} a^3  &=& \frac{1}{4}(m-1)^2m^2 \\
  \sum_{a=0}^{m-1} a^4  &=& \frac{1}{30}(m-1)m(2m-1)(3m^2-3m-1) \\
  \sum_{a=0}^{m-1} a^5  &=& \frac{1}{12}(m-1)^2 m^2(2m^2-2m-1) \\
  \sum_{a=0}^{m-1} a^6  &=& \frac{1}{42}(m-1)m(2m-1)(3m^4-6m^3+18m^2+3m+1) \\
\end{eqnarray*}
The shape of the power sum formula is an inspiration for us towards formulating the next lemma.

\bl \label{power sum divi}
Let $n$ be a positive integer.
Then we have $S_n(x) = \frac{x-1}{d_n} f_n(x)$ for a unique $f_n(x)\in\Z[x]$. Furthermore, the following statements are valid.
\begin{itemize}
  \item[$(a)$] One has $f_n(1) = d_n B_n$.
  \item[$(b)$] If $n$ is an odd integer $\geq 3$, then $x-1$ divides $f_n(x)$.
  \item[$(c)$] If $n+1$ is an odd prime, then $n+1 \nmid f_n(1)$.
\end{itemize}
\el

\bpf
 Since \[S_n(1)=\frac{1}{n+1}\big(B_{n+1}(1) - B_{n+1} \big) = \sum_{a=0}^{0} a^n = 0,\]
  the polynomial $S_n(x)$ has a factor $x-1$. Consequently, we have $S_n(x) = \frac{x-1}{d_n} f_n(x)$ for a unique $f_n(x)\in\Z[x]$. Differentiating (\ref{Bern poly}) with respect to $x$, one obtains
 \[ S_n'(x) = S_{n-1}(x) + B_n, \]
 where $S_n'(x)$ denotes the derivative of $S_n(x)$. Setting $x=1$ yields $S'_n(1) = B_n$. If $n$ is an odd integer $\geq 3$, then the latter vanishes and we obtain assertion (b). On the other hand,
 we have
 \[ S_n'(x) = \frac{f_n(x)}{d_n} + \frac{x-1}{d_n} f_n'(x). \]
 Again, setting $x=1$ and taking the above equality into account, we obtain
$f_n(1) = d_n B_n$ which is our assertion (a).

Finally, when $n+1$ is an odd prime, Proposition \ref{powersumdivisor} and Lemma \ref{Bernoulli denom} tell us that both $d_n$ and the denominator of $B_n$ are exactly divisible by $n+1$. Combining this with assertion (a), we have our assertion (c).
\epf

\br
Although we do not require this fact, we note that one can also show that $x$ divides $S_n(x)$ for every $n$ and that $x^2$ divides $S_n(x)$ for every odd $n \geq 3$. The proof, which we leave to the readers, is similar to the above.
\er

\section{Descent and Periodicity} \label{descent results}
In this section, we recall certain Galois descent property and periodicity of the even $K$-groups. As a start, we have the following descent result for extension of degree coprime to a prime $p$.

\bl \label{Kgroup}
Let $L/F$ be a finite Galois extension of number fields. Suppose that the prime $p$ does not divide the order of the Galois group $\Gal(L/F)$. Then we have
\[  K_{2k}(\Op_F)[p^\infty] \cong \Big(K_{2k}(\Op_L)[p^\infty]\Big)^{\Gal(L/F)}. \]
In particular, if $p^a$ divides $\#K_{2k}(\Op_F)$ for some $a\geq 1$, then $p^a$ also divides $\#K_{2k}(\Op_L)$.
\el

\bpf
The restriction map $K_{2k}(\Op_F)\lra K_{2k}(\Op_L)^{\Gal(L/F)}$ has kernel and cokernel annihilated by the order of the Galois group $\Gal(L/F)$ (for instance, see \cite{WeiKbook}). Since the latter order is not divisible by $p$, the conclusion of the lemma follows. Alternatively, one can identify the Sylow $p$-subgroup of $K_{2k}(\Op_F)$ with an appropriate second continuous cohomology group via the work of Rost-Voevodsky \cite{Vo} and give a cohomological proof of the said isomorphism; details of which may be found in \cite[Proposition 3.4]{LimKprank}.
\epf

We shall also require the following analog of the preceding lemma in the case of a $p$-extension.

\bl \label{Kgroup2}
Let $L/F$ be a finite Galois extension of number fields such that $|L:F|$ is a $p$-power and that $L/F$ is unramified outside $p$. Then we have
\[  K_{2k}(\Op_F)[p^\infty] \cong \Big(K_{2k}(\Op_L)[p^\infty]\Big)_{\Gal(L/F)}. \]
In particular, if $p^a$ divides $\#K_{2k}(\Op_F)$ for some $a\geq 1$, then $p^a$ also divides $\#K_{2k}(\Op_L)$.
\el

\bpf
We only give a brief sketch of the ideas of the proof. The deep work of Rost-Voevodsky \cite{Vo} allows one to identify the Sylow $p$-subgroup of $K_{2k}(\Op_F)$ with an appropriate second continuous cohomology group. Under this identification, the required isomorphism in our lemma is then a consequence of Tate spectral sequence. For more details, we refer readers to \cite[Section 2, Proposition 2.9]{KolTrieste} or \cite[Sections 3 and 4, Proposition 4.1.5]{LimKgrowth}.
\epf

In order to state the next result, we require an introduction of more terminology.
Now, if $N$ is an abelian group, we denote by $N[p]$ the subgroup of $N$ consisting of elements annihilated by $p$. Clearly, this subgroup $N[p]$ comes equipped with a $\Z/p\Z$-vector space structure, and so we may speak of the $p$-rank of $N$ which is defined by
\[ r_{p}(N) = \dim_{\Z/p\Z}\big(N[p]\big).\]
We may now state the following periodicity result.

\bp \label{periodicity of p-rank}
Let $F$ be a number field. Then we have
\[ r_p\big( K_{2k}(\Op_F)\big) = r_p\big( K_{2k'}(\Op_F)\big), \]
whenever $k \equiv k'$ $($mod $|F(\ze_p):F|)$. In particular, one has $p \mid \#K_{2k}(\Op_F)$
if and only if $p \mid \#K_{2k'}(\Op_F)$.
\ep

\bpf
This is well-known among the experts. For the convenience of the readers, we sketch a proof here. Thanks to the deep work of Rost and Voevodsky \cite{Vo}, there is an identification
\[ K_{2k}(\Op_F)/p \cong H^2\big(\Gal(F_{S_p}/F)), \mu_p^{\otimes(k+1)}\big),\]
where $F_{S_p}$ is the maximal algebraic extension of $F$ unramified outside the set of primes of $F$ above $p$. If $j \equiv 0$ mod $[F(\zeta_p) : F]$, then the Galois group $\Gal(F_{S_p}/F)$ acts trivially on $\mu_p^{\otimes j}$. Therefore, it follows that if $k \equiv
k'$ mod $[F(\zeta_p) : F]$, then we have
\[
  \begin{split}
    K_{2k}(\Op_F)/p  & \cong H^2\big(\Gal(F_{S_p}/F)), \mu_p^{\otimes(k+1)}\big) \\
       &\cong H^2\big(\Gal(F_{S_p}/F)), \mu_p^{\otimes(k'+1)}\big)\otimes \mu_p^{\otimes(k-k')}\\
       &\cong K_{2k'}(\Op_F)/p\otimes \mu_p^{\otimes(k-k')}.
   \end{split}\]
  Consequently, the groups $K_{2k}(\Op_F)/p$ and $K_{2k'}(\Op_F)/p$ have the same rank over $\Z/p\Z$.
\epf

\section{Lichtenbaum conjecture for totally real abelian number fields} \label{Lich sec}

From now on, we always let $\ze_F(s)$ denote the Dedekind zeta function of a number field $F$. It is well-known that this said function has an analytical continuation to the whole complex plane except at the point $s=1$ at which it has a simple pole. Therefore, it makes sense to speak of $\zeta_F(-k)$ for a positive integer $k$. Let $\mu_\infty$ be the group of all the roots of unity of $\bar{F}$, where $\bar{F}$ is the algebraic closure of $F$. For an integer $j\geq 1$, we write $\mu_\infty^{\otimes j}$ for the $j$-fold tensor product of $\mu_\infty$ with $\Gal(\bar{F}/F)$ acting diagonally. Set $w_j(F)$ to be the order of $(\mu_\infty^{\otimes j})^{\Gal(\bar{F}/F)}$ for a number field $F$.

We can now state the following deep result which will be required in the discussion of the paper.

\bt \label{Lich}
Let $F$ be a totally real abelian number field of degree $r$. Then for every odd integer $k\geq 1$, one has
 \[ \#K_{2k}(\Op_F) =\begin{cases} (-1)^{r}w_{k+1}(F)\zeta_F(-k),  & \mbox{if $k\equiv 1$ (mod 4)}, \vspace{.1in}\\
 \displaystyle\frac{1}{2^r}w_{k+1}(F)\zeta_F(-k), & \mbox{if $k\equiv 3$ (mod 4)}.
 \end{cases} \]
\et

\bpf
This identity was first conjectured by Birch and Tate for $K_2$ (see \cite{Tate}) and subsequently, Lichtenbaum extended their conjecture to general $K_{2k}$ \cite{Lic72, Lic73}. It was probably due to the insight of Coates \cite{C72, C73} that one might possibly prove this conjecture via the main conjecture of Iwasawa \cite{Iw73}. Later, B\'{a}yer and Neukirch \cite{BN} confirmed this vision of Coates by showing that the main conjecture of Iwasawa implies a cohomological version of Lichtenbaum's conjecture (we refer readers to \cite{BN} for the precise formulation of this cohomological Lichtenbaum's conjecture). This implication result of B\'{a}yer and Neukirch is concerned with odd primes. The implication for $p=2$ is due to Kolster (see \cite{Kol}; also see appendix in \cite{RW}). The cohomological version of Lichtenbaum's conjecture is equivalent to the $K$-theoretical version via the Quillen-Lichtenbaum conjecture (see \cite{KolTrieste}), which is now a theorem being a consequence of the monumental work of Rost-Voevodsky (\cite{Vo}; also see the previous  important works of Soul\'e \cite{Sou}, Dwyer-Friedlander \cite{DF} and Rognes-Weibel \cite{RW} related to this). Prior to this, the main conjecture of Iwasawa has been proved by Mazur-Wiles \cite{MW} and Wiles \cite{Wiles} (also see \cite{Ban, Ban2, NQD}).
\epf

\section{Proof of the main result} \label{Main results section}

This section is concerned with proving the main results of the paper. As a start, let $\chi$ be a nontrivial even Dirichlet character of order $p^b$ and conductor $m$. Recall that the
generalized Bernoulli numbers $B_{n,\chi}$ are defined by
\[ \sum_{a=1}^m\frac{\chi(a)te^{at}}{e^{mt}-1} = \sum_{n=0}^\infty B_{n,\chi}\frac{t^n}{n!}.\]
It follows from \cite[Proposition 4.1]{Wash} that one has
\[ B_{k+1,\chi} = m^{k}\sum_{a=1}^{m}\chi(a)B_{k+1}(a/m).  \]
Since the order of $\chi$ is $p^b$, the values of $\chi$ are $p^b$-th roots of unity. Therefore, it follows that $B_{k+1,\chi}$ may be viewed as an element in $\Q(\ze_{p^N})$ for every $N\geq b$.

For a nonzero element $\al \in \Q(\ze_{p^N})$, we may write $\al =\be /r$ for some element $\be\in \Z[\ze_{p^N}]$ and integer $r$. In the event that $r$ is coprime to $p$, it then makes sense to speak of $\al$ $\mathrm{mod}~ (1-\ze_{p^N})\Z[\ze_{p^N}]$ which is essentially
\[ \be r^{-1} ~\mathrm{mod}~ (1-\ze_{p^N})\Z[\ze_{p^N}].\]
In particular, if $\be \in (1-\ze_{p^N})\Z[\ze_{p^N}]$ and $r$ is coprime to $p$, we shall write
\[ \al \equiv 0\quad\mathrm{mod}~ (1-\ze_{p^N})\Z[\ze_{p^N}]. \]
This convention will be consistently applied throughout the subsequent discussion without further mention. Equipped with this understanding, we are now in a position to present the following result, which will serve as a crucial component for the eventual proof of our main theorem.

\bp \label{mainchar1} Let $\chi$ be a nontrivial even Dirichlet character of order $p^b$ and conductor $m$, where $p$ is an odd prime and $m$ is a squarefree integer $>1$ with the property that every of its prime divisor satisfies $\equiv 1$ $(\mathrm{mod}$ $p)$.
Suppose that either one of the following statements is valid.
\begin{itemize}
\item[$(i)$] $p>k+2$.
\item[$(ii)$] $p=k+2$ and at least one prime divisor $\ell$ of $m$ satisfies $v_p(l-1)\geq 2$.
\end{itemize}
Then we have $B_{k+1,\chi}\equiv 0 ~ \mathrm{mod}~(1-\ze_{p^N})^{p^{N-b}}\Z[\ze_{p^N}]$ for every $N\geq b$.
\ep

\bpf
 To begin, we set $m = \ell_1\cdots \ell_s$, where the $\ell_i$'s are distinct primes satisfying $p\mid \ell_i-1$.
Since $(1-\ze_{p^b})\Z[\ze_{p^N}] = (1-\ze_{p^N})^{p^{N-b}}\Z[\ze_{p^N}]$, it suffices to prove the congruence for $N=b$. By \cite[Proposition 4.1]{Wash}, one has
\[ B_{k+1,\chi} = m^{k}\sum_{a=1}^{m-1}\chi(a)B_{k+1}(a/m)  \]
 noting that $\chi(m)=0$. The sum on the right can then be rewritten as
 \begin{eqnarray*}
m^{k}\sum_{a=1}^{m-1}\chi(a)\left(\sum_{i=0}^{k+1} \binom{k+1}{i}(a/m)^{k+1-i}B_i \right)
 &=& \sum_{i=0}^{k-1}\binom{k+1}{i}B_i m^{i-1} \left(\sum_{a=1}^{m-1} \chi(a)a^{k+1-i}\right),
                                                                  \end{eqnarray*}
 where we have made use of the fact that the sums $\sum_{a=1}^{m-1} \chi(a)$ and $\sum_{a=1}^{m-1} \chi(a)a$ both vanish, due to the fact that $\chi$ is a nontrivial even Dirichlet character. As the order of $\chi$ is $p^b$, the values of $\chi$ are $p^b$-th roots of unity. Therefore, the sum $\sum_{a=1}^{m-1}\chi(a)a^{k+1-i}$ can be viewed as an element in the ring $\Z[\ze_{p^b}]$.  Consequently, we have
\begin{equation}\label{B=powersum}
   B_{k+1,\chi} \equiv \sum_{i=0}^{k-1}\binom{k+1}{i}B_i m^{i-1} \left(\sum_{\substack{a=1 \\ \gcd(a,m)=1}}^{m-1} a^{k+1-i}\right) \quad\mathrm{mod}~ (1-\ze_{p^b})\Z[\ze_{p^b}]  \end{equation}
noting that the congruence makes sense in view of our assumption that $p\geq k+2$ and that the Bernoulli numbers appearing have denominators with prime divisors $\leq k$ (cf. Lemma \ref{Bernoulli denom}). Now observe that
\begin{equation}\label{B=powersum2} \sum_{\substack{a=1 \\ \gcd(a,\ell_1\cdots \ell_s)=1}}^{\ell_1\cdots \ell_s-1} a^{k+1-i} = \sum_{a=1 }^{\ell_1\cdots \ell_s-1} a^{k+1-i} + \sum_{t=1}^{s-1}(-1)^t\left(\sum_{1\leq j_1<\cdots< j_t\leq s}(\ell_{j_1}\cdots \ell_{j_t})^{k+1-i}\sum_{a=1}^{\frac{\ell_1\cdots \ell_s}{\ell_{j_1}\cdots \ell_{j_t}}-1} a^{k+1-i}\right). \end{equation}
In view of (\ref{B=powersum}) and (\ref{B=powersum2}), we are essentially reduced to showing that sum of the form
\begin{equation}\label{powersum of interest}\sum_{a=1 }^{\ell_1\cdots \ell_s-1} a^{k+1-i}
\end{equation}
is divisible by $p$. By Lemma \ref{power sum divi}, the latter sum is precisely
\[ (\ell_1\cdots \ell_s-1) \frac{f_{k+1-i}(\ell_1\cdots \ell_s)}{d_{k+1-i}}. \]
 We first consider the case that $p>k+2$. It then follows from Lemma \ref{power sum divi}(c) that $d_{k+1-i}$ is coprime to $p$ for $i=0,1,...,k-2$. On the other hand, since each $\ell_j-1$ is divisible by $p$, so is $\ell_1\cdots \ell_s-1$. This therefore proves that the sum (\ref{powersum of interest}) is divisible by $p$ when $p>k+2$.

We now turn to the case when $p=k+2$.
The argument in the preceding paragraph shows that the sum (\ref{B=powersum2}) is divisible by $p$ for $i=1,...,k-1$. Therefore, it remains to verify that
\[\sum_{\substack{a=1 \\ \gcd(a,\ell_1\cdots \ell_s)=1}}^{\ell_1\cdots \ell_s-1} a^{k+1}\]
is divisible by $p$. Since $p \parallel d_{k+1}$ by Proposition \ref{powersumdivisor}(c), we are reduced to showing that
\begin{equation}\label{sum p^2}
  d_{k+1}\sum_{\substack{a=1 \\ \gcd(a,\ell_1\cdots \ell_s)=1}}^{\ell_1\cdots \ell_s-1} a^{k+1}
\end{equation}
is divisible by $p^2$. Relabelling, if necessary, we may assume that $\ell_s\equiv 1$ (mod $p^2$). We then proceed as in (\ref{B=powersum2}) and apply Lemma \ref{power sum divi} to obtain
\begin{eqnarray*}
 d_{k+1}\sum_{\substack{a=1 \\ \gcd(a,\ell_1\cdots \ell_s)=1}}^{\ell_1\cdots \ell_s-1} a^{k+1}  &=& (\ell_1\cdots \ell_s-1)f_{k+1}(\ell_1\cdots \ell_s) - \ell_s^{k+1} (\ell_1\cdots \ell_{s-1}-1)f_{k+1}(\ell_1\cdots \ell_{s-1}) \\
   &+&  \sum_{t=1}^{s-2}(-1)^t\left(\sum_{1\leq j_1<\cdots< j_t\leq s-1}(\ell_{j_1}\cdots \ell_{j_t})^{k+1}\Big(\frac{\ell_1\cdots \ell_s}{\ell_{j_1}\cdots \ell_{j_t}}-1 \Big)f_{k+1}\Big(\frac{\ell_1\cdots \ell_s}{\ell_{j_1}\cdots \ell_{j_t}}\Big)\right)\\
 &+&  \sum_{t=1}^{s-2}(-1)^{t+1}\left(\sum_{1\leq j_1<\cdots< j_t\leq s}(\ell_{j_1}\cdots \ell_{j_t}\ell_s)^{k+1}\Big(\frac{\ell_1\cdots \ell_s}{\ell_{j_1}\cdots \ell_{j_t}\ell_s}-1 \Big)f_{k+1}\Big(\frac{\ell_1\cdots \ell_s}{\ell_{j_1}\cdots \ell_{j_t}\ell_s}\Big)\right)\\
 & & \\
 &+& (-1)^{s-1}(\ell_1\cdots \ell_{s-1})^{k+1}(\ell_s-1)f_{k+1}(\ell_s).
    \end{eqnarray*}
Using the fact that $\ell_s\equiv 1$ (mod $p^2$), we see that the two terms in the first line cancel off modulo $p^2$. Similarly, one can check that the sum in the second line cancels off the sum in the third line modulo $p^2$. The final term in the last line is plainly $\equiv 0$ (mod $p^2$). Thus, we conclude that the expression (\ref{sum p^2}) is divisible by $p^2$. This completes the proof of the proposition.
\epf

We are in position to prove the following theorem.

\bt \label{main} Let $p$ be an odd prime, $k$ a positive odd integer and $m$ an integer $>1$. Define
\[s_j = \#\big\{\ell \mbox{ is a prime divisor of } m~\big|~ v_p(\ell-1)=j\big\} \]
and set $\theta:=\max\{j~|~s_j\neq 0\}$. Suppose that $p\geq k+2$.
Then the power of $p$ dividing the orders of the groups $K_{2k}(\Z[\ze_m+\ze_m^{-1}])$ is at least
\[\frac{p^{s_1+s_2\cdots +s_\theta}-p^{\delta s_1}+\delta}{p-1} + \sum_{j=2}^{\theta}\frac{p^{s_j+\cdots +s_\theta}-1}{p-1}p^{s_1+2s_2+\cdots + (j-1)s_{j-1} + (j-1)(s_j+\cdots s_\theta) -(j-1)},\]
where
\[ \delta =\left\{
     \begin{array}{ll}
       0, & \hbox{if $p>k+2$;} \\
       1, & \hbox{if $p=k+2$.}
     \end{array}
   \right.\]
A similar divisibility assertion also holds for $K_{2k}(\Z[\ze_m])$.
\et

\bpf
Since the prime $p$ is odd, Lemma \ref{Kgroup} tells us that it suffices to prove the divisibility assertion for $K_{2k}(\Z[\ze_m+\ze_m^{-1}])$. Now suppose that $p$ divides $m$. We then write $m = p^am'$, where $a\geq 1$ and $p\nmid m'$. Note that $\Q(\ze_m)/\Q(\ze_{pm'})$ is a finite $p$-extension which is unramified outside the primes above $p$. On the other hand, $\Q(\ze_{pm'})/\Q(\ze_{m'})$ is a Galois extension of degree coprime to $p$. Therefore, it follows from Lemmas \ref{Kgroup} and \ref{Kgroup2} that for the verification of the divisibility assertion for $K_{2k}(\Z[\ze_m+\ze_m^{-1}])$, it is enough to verify the divisibility assertion for $K_{2k}(\Z[\ze_{m'}+\ze_{m'}^{-1}])$.

Hence, by virtue of the discussion in the preceding paragraph, we may, and will, assume that the integer $m$ we start with is not divisible by the prime $p$. Taking Lemma \ref{Kgroup} into account again, we may even assume that $m$ is squarefree and that every prime divisor $\ell$ of $m$ satisfies $\ell\equiv 1$ mod $p$. Let $F$ be the maximal $p$-extension of $\Q$ contained in $\Q(\ze_{m}+\ze_{m}^{-1})$. Via Lemma \ref{Kgroup} again, we are reduced to prove the divisibility assertion for $K_{2k}(\Op_F)$.  Note that the values of $s_j$ are unchanged throughout the whole reduction process.  By our assumption on $m$, we see that
\[ G:=\Gal(F/\Q)\cong \bigoplus_{\ell\mid m}\Z/p^{v_p(\ell-1)}\Z =\bigoplus_{j=1}^{\theta} \big(\Z/p^j\Z\big)^{\oplus s_j}, \]
and so the group $\widehat{G}$ of linear characters of $G$ has the same group structure. Set $\# \widehat{G}=p^N$. In particular, we have
\[ N = s_1+2s_2+\cdots+ \theta s_{\theta}.\]

As $F$ is a totally real abelian field, we may apply Theorem \ref{Lich} to see that
\[ \#K_{2k}(\Op_F) \sim_p w_{k+1}(F)\zeta_F(-k).\]
By Artin formalism, there is the following factorization
\[ \ze_F(-k) = \ze(-k) \prod_{\substack{\chi\in \widehat{G}\\ \chi\neq\chi_0}}L(\chi, -k), \]
where $\chi_0$ denotes the trivial character and $L(\chi, s)$ is the Dirichlet $L$-function attached to $\chi$. It's well-known that
\[ \ze(-k)=-\frac{B_{k+1}}{k+1} \quad\mathrm{and}\quad L(\chi, -k) = -\frac{B_{k+1,\chi}}{k+1} \]
(for instance, see \cite[Theorem 4.2]{Wash}). Therefore, we have
\begin{equation}\label{K=Bernoulli}
  \# K_{2k}(\Op_F) \sim_p \frac{w_{k+1}(F) B_{k+1}}{(k+1)^{p^N}} \prod_{\substack{\chi\in \widehat{G}\\ \chi\neq\chi_0}}B_{k+1,\chi}.
\end{equation}

We claim that
\begin{equation} \label{B=p 1}
  \frac{w_{k+1}(F) B_{k+1}}{(k+1)^{p^N}}\sim_p 1.
\end{equation}
Indeed, since $p\geq k+2$, the denominator term has no $p$ contribution. On the other hand, it follows from \cite[Chap.\ VI, Proposition 2.2]{WeiKbook} that the prime $p$ divides $w_{k+1}(F)$ if and only $p-1 \mid k+1$. Furthermore, the same proposition tells us that $p\parallel w_{k+1}(F)$ in the context of $p=k+2$. Hence we conclude that $p\nmid w_{k+1}(F)$ or $p\parallel w_{k+1}(F)$ accordingly to $p>k+2$ or $p=k+2$. Since this $p$-divisibility property also holds for the denominator of $B_{k+1}$ by Lemma \ref{Bernoulli denom}, our claim (\ref{B=p 1}) thus follows.

Now, combining (\ref{K=Bernoulli}) and (\ref{B=p 1}), we obtain
\begin{equation}\label{K=Bernoulli2}
  \# K_{2k}(\Op_F) \sim_p  \prod_{\substack{\chi\in \widehat{G}\\ \chi\neq\chi_0}}B_{k+1,\chi}.
\end{equation}
Set $\widehat{G}_j:=\{\chi\in \widehat{G}~|~\mathrm{ord}(\chi)=p^j\}$. Then we have the following factorization
\begin{equation} \label{B fact} \prod_{\substack{\chi\in \widehat{G}\\ \chi\neq\chi_0}}B_{k+1,\chi} =  \prod_{j=1}^\theta\Big( \prod_{\chi\in \widehat{G}_j}B_{k+1,\chi}\Big).
\end{equation}
It follows from Proposition \ref{mainchar1} that
\begin{equation} \label{B prod=p}\prod_{\chi\in \widehat{G}_j}B_{k+1,\chi}  \equiv 0~\mathrm{mod}~ (1-\ze_{p^N})^{p^{N-j}\#\widehat{G}_j}\Z[\ze_{p^N}]\end{equation}
for $j\geq 2$. When $j=1$, the same proposition also tells us that
\begin{equation} \label{B prod=p1}\prod_{\chi\in \widehat{G}_1}B_{k+1,\chi}  \equiv 0~\mathrm{mod}~ (1-\ze_{p^N})^{p^{N-1}\#\widehat{G}_1}\Z[\ze_{p^N}]\end{equation}
provided that $p>k+2$.

Now suppose that $p=k+2$. Write $m = m^+m^-$, where every prime divisor of $m^+$ satisfies $\equiv 1$ (mod $p^2$). By Proposition \ref{mainchar1}, we see that $B_{k+1,\chi}  \equiv 0~\mathrm{mod}~ (1-\ze_{p^N})^{p^{N-1}}\Z[\ze_{p^N}]$ whenever the conductor of $\chi$ contains a prime divisor $\equiv 1$ (mod $p^2$). On the other hand, characters with conductors not having this property must factor through $\Gal(\Q(\ze_{m^-})/\Q)$, and so there are exactly $p^{s_1}-1$ of these. Also, note that since such characters must necessarily have order $p$, we therefore only need to count the contribution coming from these characters in $\widehat{G}_1$, and so we see that
\begin{equation} \label{B prod=p2} \prod_{\chi\in \widehat{G}_1}B_{k+1,\chi}  \equiv 0~\mathrm{mod}~ (1-\ze_{p^N})^{(\#\widehat{G}_1-p^{s_1}+1)p^{N-1}}\Z[\ze_{p^N}]\end{equation}

Combining (\ref{K=Bernoulli2})$-$(\ref{B prod=p2}), we have
\[ \# K_{2k}(\Op_F)\in (1-\ze_{p^N})^{p^{N-1}(\#\widehat{G}_1-p^{\delta s_1}+\delta) +\sum_{j=2}^\theta p^{N-j}\#\widehat{G}_j}\Z[\ze_{p^N}]\cap \Z = p^w\Z, \]
where
\begin{equation} \label{w=sum} w = \left\lceil \frac{p^{N-1}(\#\widehat{G}_1-p^{\delta s_1}+\delta) +\sum_{j=2}^\theta p^{N-j}\#\widehat{G}_j}{p^{N-1}(p-1)}\right\rceil  = \left\lceil \frac{\#\widehat{G}_1-p^{\delta s_1}+\delta}{p-1} + \sum_{j=2}^\theta \frac{\#\widehat{G}_j}{p^{j-1}(p-1)}\right\rceil. \end{equation}
By an elementary abelian $p$-group calculation, one has
\[\#\widehat{G}_j = p^{s_1+ 2s_2+\cdots +(j-1)s_{j-1}+ j(s_j+ \cdots + s_\theta)} - p^{s_1+ 2s_2+\cdots +(j-1)(s_{j-1} + s_j+ \cdots + s_\theta)}. \]
Putting this into (\ref{w=sum}), it then follows from a direct calculation that we obtain the quantity as asserted in the theorem.
\epf

\br \label{remark div}
We discuss two situations, where the expression in Theorem \ref{main} simplifies.
\begin{enumerate}
  \item[$(1)$] Suppose that $\theta =1$ and $p >k+2$. Then in this situation, the power of $p$ in Theorem \ref{main} has the simple form
\[ \frac{p^{v_p(m)}-1}{p-1}. \]
For instances, one has
\[ 5^{1+5+5^2+5^3+5^4} \mid \#K_{2}(\Z[\zeta_{11\cdot 31\cdot 41 \cdot 61 \cdot 71 }]) \quad\mbox{and}\quad 7^{1+7+7^2} \mid \#K_{6}(\Z[\zeta_{29\cdot 43\cdot 71}]). \]

  \item[$(2)$] Suppose that $\theta \geq 2$ and $s_j =0$ for $j\neq \theta$. Then the power of $p$ is
 \[\frac{p^{s_\theta}-1 +\delta}{p-1} + \sum_{j=2}^{\theta}\frac{p^{s_\theta}-1}{p-1}p^{(j-1) (s_\theta -1)},\]
\end{enumerate}
\er

We conclude with the following.

\bpf[Proof of Theorem \ref{mainintro}]
The conclusion of the theorem follows from a combination of Theorem \ref{main} and Proposition \ref{periodicity of p-rank}.
\epf

\section{Revisiting and generalizing a result of Browkin} \label{BrowkinSec}

In this section, we return to study the situation of $p=k+2$ in greater depth for a Dirichlet character of prime conductor. Our main result is as follows.

\bp \label{Browkinchar} Let $\chi$ be a nontrivial even Dirichlet character of order $p$ and conductor $\ell$, where $p$ is an odd prime and $\ell$ is a prime with the property that $p \mid \ell-1$.
Then $B_{p-1,\chi}\equiv 0 ~ \mathrm{mod}~(1-\ze_{p})\Z[\ze_{p}]$ if and only if $v_p(\ell-1)\geq 2$.
\ep

\bpf
Write $k=p-2$.
By a similar argument to that in Proposition \ref{mainchar1}, we have
\begin{eqnarray*}
 B_{k+1,\chi}  & \equiv & \sum_{i=0}^{k-1}\binom{k+1}{i}B_i \ell^{i-1} \left(\sum_{a=1}^{\ell-1} a^{k+ \ell -i}\right) \quad\mathrm{mod}~ (1-\ze_{p})\Z[\ze_{p}] \\
 &\equiv & \sum_{i=0}^{k-1}\binom{k+1}{i}B_i \ell^{i-1} \frac{(\ell-1) f_{k+ \ell -i}(\ell )}{\ell d_{k+ 1 -i}} \quad\mathrm{mod}~ (1-\ze_{p})\Z[\ze_{p}]
\end{eqnarray*}
The sum on the right can be rewritten as
\[ \frac{(\ell-1) f_{k+1}(\ell)}{\ell d_{k+1}} + \sum_{i=1}^{k-1}\binom{k+1}{i}B_i \ell^{i-1} \frac{(\ell-1) f_{k+\ell-i}(\ell)}{\ell d_{k+1-i}}.\]
As seen in the proof of Proposition \ref{mainchar1}, the sum on the right is divisible by $p$. Hence it follows that  $B_{k+1,\chi}\equiv 0 ~ \mathrm{mod}~(1-\ze_{p})\Z[\ze_{p}]$ if and only if $\frac{(l-1) f_{k+1}(l)}{ld_{k+1}}$ is divisible by $p$. Now, Proposition \ref{powersumdivisor}(c) tells us that $p\parallel d_{k+1}$. Furthermore, it follows from Lemma \ref{power sum divi}(c) that
\[ f_{k+1}(l)\equiv f_{k+1}(1)  \not\equiv 0~(\mathrm{mod}~p). \]
 Hence we see that  $\frac{(l-1) f_{k+1}(l)}{ld_{k+1}}$ is divisible by $p$ if and only if  $p^2 \mid l-1$. This completes the proof of the proposition.
\epf

We translate the preceding result to the context of $K$-groups.

\bt \label{BrowkinK} Let $F$ be a cyclic extension of $\Q$ with $|F:\Q|=p$ and $F\subseteq \Q(\ze_\ell)$ for some prime $\ell$.
Then $p$ divides $\#K_{2(p-2)}(\Op_F)$ if and only if $v_p(\ell-1)\geq 2$.
\et

\bpf By (\ref{K=Bernoulli2}), we have
\begin{equation}\label{K=Bernoulli3}
  \# K_{2k}(\Op_F) \sim_p \prod_{\substack{\chi\in \widehat{\Gal(F/\Q)}\\ \chi\neq\chi_0}}B_{k+1,\chi}.
\end{equation}
Therefore, the prime $p$ divides $\#K_{2k}(\Op_F)$ if and only if $B_{k+1,\chi}\equiv 0 ~ \mathrm{mod}~(1-\ze_{p^b})\Z[\ze_{p^b}]$ for some $\chi\in \widehat{\Gal(F/\Q)}-\{\chi_0\}$. The latter is precisely equivalent to $v_p(\ell-1)\geq 2$ by Proposition \ref{Browkinchar}.
\epf

\br
  When $p=3$, Theorem \ref{BrowkinK} was proved by Browkin \cite[Theorem 2.4]{Brow05}. Therefore, our theorem can be thought as a generalization of Browkin's result to higher $K$-groups. Proposition \ref{Browkinchar} and Theorem \ref{BrowkinK} have also been proved in \cite[Section 6.1]{AssMov} by a different approach (we thank the anonymous referee for pointing this out). 
\er

We record the following quantitative observation which gives us a sense of how often the divisibility property in Theorem \ref{BrowkinK} is valid. We should mention that the corollary can also be deduced from the discussion in \cite[Corollary 4.5 and Section 6.1]{AssMov}.

\bc \label{BrowkinKcor}
Let $p$ be an odd prime. For every prime $\ell\equiv 1$ $($mod $p)$, we write $F_\ell$ for the  cyclic extension of $\Q$ of degree $p$ contained in $\Q(\ze_\ell)$. Then one has
\[\lim_{x\to\infty}\frac{\#\{\ell~prime \leq x~|~\mbox{$p$ divides $\#K_{2(p-2)}(\Op_{F_\ell})$ and $\ell\equiv 1$ $($mod $p)$}\}}{\#\{\ell~prime \leq x~|~\mbox{$\ell\equiv 1$ $($mod $p)$}\}} =\frac{~1~}{~p~}.\]
\ec

\bpf
Let $\pi_a(x) = \#\{\ell~\mathrm{prime} \leq x~|~l\equiv 1 ~(\mathrm{mod}~ a)\}$ for $a=p, p^2$. By Proposition \ref{BrowkinK}, we have
\[ \pi_{p^2}(x) =  \#\{\ell~\mathrm{prime} \leq x~|~\mbox{$p$ divides $\#K_{2k}(\Op_{F_\ell})$ and $\ell\equiv 1$ $($mod $p)$}\}.\]
On the other hand, the theorem of Dirichlet arithmetic progression tells us that
\[ \lim_{x\to\infty}\frac{\pi_a(x)}{\pi(x)} =\frac{1}{\phi(a)}\]
for $a=p, p^2$, where $\pi(x)$ is the number of primes $\leq x$ and $\phi$ is the Euler totient function. Therefore, it then follows that
\[ \lim_{x\to\infty}\frac{\pi_{p^2}(x)}{\pi_p(x)} = \lim_{x\to\infty}\frac{\pi_{p^2}(x)}{\pi(x)} \frac{\pi(x)}{\pi_p(x)} = \frac{\phi(p)}{\phi(p^2)} =\frac{~1~}{~p~} \]
as required.
\epf

\section{Examples} \label{example sec}

In this final section, we give some examples to illustrate the results of this paper. All the computations of $K$-groups in this section is done appealing to Theorem \ref{Lich} with the aid of Pari/GP \cite{Pari}.

\begin{itemize}
\item By Proposition \ref{BrowkinK}, $\#K_2(\Z[\ze_{7}+\ze_{7}^{-1}])$ is not divisible by $3$. It then follows from Proposition \ref{periodicity of p-rank} that $\#K_{2k}(\Z[\ze_{7}+\ze_{7}^{-1}])$ is not divisible by $3$ for every odd positive integer $k$. Via Pari, we have

$\#K_2(\Z[\ze_{7}+\ze_{7}^{-1}])=8=2^3$,

$\#K_6(\Z[\ze_{7}+\ze_{7}^{-1}])=79$,

 $\#K_{10}(\Z[\ze_{7}+\ze_{7}^{-1}])=59144=2^3\cdot 7393$,

$\#K_{14}(\Z[\ze_{7}+\ze_{7}^{-1}])=142490119$.

$\#K_{18}(\Z[\ze_{7}+\ze_{7}^{-1}])=913161859868= 2^2\cdot 228290464967$.

$\#K_{22}(\Z[\ze_{7}+\ze_{7}^{-1}])= 2101941875088322867= 691\cdot 10903\cdot 278995143079$.

As seen from these computations, these $K$-groups have orders not divisible by $3$.

\item Proposition \ref{BrowkinK} tells us that $\#K_{8t+2}(\Z[\ze_{11}+\ze_{11}^{-1}])$ is divisble by $5$, whereas $\#K_{8t+6}(\Z[\ze_{11}+\ze_{11}^{-1}])$ is not. Via Pari, we computed a few $K$-groups and see that they do satisfy this expected property.

$\#K_2(\Z[\ze_{11}+\ze_{11}^{-1}])=160=2^5\cdot 5$,

$\#K_6(\Z[\ze_{11}+\ze_{11}^{-1}])=847811=71\cdot 11941$,

$\#K_{10}(\Z[\ze_{11}+\ze_{11}^{-1}])=407495402731360=2^5\cdot 5\cdot 521\cdot 4888380551$,

$\#K_{14}(\Z[\ze_{11}+\ze_{11}^{-1}])=3543010400763352360091= 13721\cdot 2520121\cdot  102462575851$.


\item It follows from Propositions \ref{periodicity of p-rank}  and \ref{BrowkinK} that $\#K_{2k}(\Z[\ze_{13}+\ze_{13}^{-1}])$ is not divisible by $3$ for every odd positive integer $k$. The following computations confirm this.

$\#K_2(\Z[\ze_{13}+\ze_{13}^{-1}])=1216=2^6\cdot 19$,

$\#K_6(\Z[\ze_{13}+\ze_{13}^{-1}])=316792259 = 7\cdot 29\cdot 103 \cdot 109\cdot 139$

$\#K_{10}(\Z[\ze_{13}+\ze_{13}^{-1}])= 99222088525421989696 = 2^6 \cdot 73\cdot 109\cdot 307\cdot 2341 \cdot 2953 \cdot 91807$

  \item On the other hand, Propositions \ref{periodicity of p-rank}  and \ref{BrowkinK} tell us that $\#K_{2k}(\Z[\ze_{19}+\ze_{19}^{-1}])$ is divisible by $3$ for every odd positive integer $k$. The following numerical computations illustrate this.

 $\#K_2(\Z[\ze_{19}+\ze_{19}^{-1}])=2244096 = 2^9\cdot 3^2\cdot 487$.

$\#K_{6}(\Z[\ze_{19}+\ze_{19}^{-1}])= 540700931767472649 = 3^2\cdot 61 \cdot 67 \cdot 883 \cdot 16647509341$




\item The next pair of numerical examples illustrates the $11$-divisibility of $K$-groups of $\Z[\ze_{23}+\ze_{23}^{-1}]$.

$\#K_2(\Z[\ze_{23}+\ze_{23}^{-1}])=837613568 =2^{11}\cdot 11 \cdot 37181$,

$\#K_6(\Z[\ze_{23}+\ze_{23}^{-1}])= 6952891386341432645005057 = 11\cdot 1607 \cdot 120263419 \cdot 3270569157439$


\item We computed

$\#K_2(\Z[\ze_{31}+\ze_{31}^{-1}])=580922038681600=2^{17}\cdot 5^2\cdot 7 \cdot 11 \cdot 2302381$

which is in line with our expectation that its order is divisible by $5$.
    \end{itemize}

\end{document}